\documentclass[12pt]{article}

\usepackage{amsfonts}
\usepackage{amssymb}
\usepackage{amsmath}
\usepackage{amsthm}

\def \le {\leqslant}
\def \ge {\geqslant}

\begin{document}

\begin{Large}
\centerline{\bf Proof of W.M.Schmidt's conjecture }
\centerline{\bf concerning successive minima of a lattice}
\vskip+1.5cm \centerline{\bf Moshchevitin N.G. \footnote{ The
research was supported by the grant RFBR 09-01-00371a
 }
}
\end{Large}
\vskip+1.5cm

\centerline{\bf Abstract} We prove W.M.Schmidt's conjecture  about
a one-parameter family of lattices related to simultaneous
Diophantine approximations.

 \vskip+1.5cm

\section{Introduction}
  Consider real numbers $\xi_j \in[0,1), 1\le j\le n$. For a real   $x$ we denote by $|x|$  the absolute value of $x$.
We know since Dirichlet that for any real $N>1$ the inequalities
$$
|x|\le N,\,\,\, \max_{1\le j\le n}|x\xi_j -y_j|\le N^{-1/n}
$$
have a solution in integers $x\neq 0, y_1,...,y_n$. Similarly for
any real $N>1$ the inequalities
$$
|x|\le N,\,\,\, \left(\sum_{j=1}^n|x\xi_j -y_j|^2\right)^{1/2}\le
2w_n^{-1/n} N^{-1/n}
$$
have a solution in integers $x\neq 0, y_1,...,y_n$ (here  $w_n$
stands for  the volume of the unit ball in the $n$-dimensional
Euclidean space).

In this paper we work with Euclidean space $\mathbb{R}^{n+1}$ with
coordinates $(x,y_1,...,y_n)$ and with Euclidean space
$\mathbb{R}^{n}$ with coordinates $(y_1,...,y_n)$.

Consider an $(n+1)$-dimensional vector $ \xi = (1,\xi_1,...,\xi_n)
\in \mathbb{R}^{n+1}$.

For a vector $y =(y_1,...,y_n) \in \mathbb{R}^n$ we define $|y | $
to be the Euclidean norm of $y$. So, $|y|
=\sqrt{y_1^2+...+y_n^2}.$  We also use the notation $|y|_s =
\max_{1\le j\le n} |y_j|$ for the sup-norm of a vector $y
=(y_1,...,y_n) \in \mathbb{R}^n$.

 For a
real $N\ge 1$ and a vector $\xi $  we define a matrix
$$
{\cal A} (\xi, N) = \left(
\begin{array}{ccccc}
N^{-1} & 0& 0&  \cdots &0 \cr N^{\frac{1}{n}} \xi_1 &
-N^{\frac{1}{n}} & 0&\cdots & 0 \cr  N^{\frac{1}{n}} \xi_2 &0&
-N^{\frac{1}{n}} &  \cdots & 0 \cr \cdots &\cdots &\cdots &\cdots
\cr N^{\frac{1}{n}} \xi_n &0&0&\cdots &- N^{\frac{1}{n}}
\end{array}\right)
$$
and a lattice
$$
\Lambda (\xi, N) = {\cal A} (\xi, N)\mathbb{Z}^{n+1}.
$$
Consider the $ (n+1)$-dimensional unit cube
$${\cal U} =
\{ z= (x,y_1,...,y_n)\in \mathbb{R}^{n+1}:\,\,\, \max ( |x|, |y|_s
)\le 1  \}
$$
and a convex 0-symmetric body
$${\cal W} =
\{ z= (x,y_1,...,y_n)\in \mathbb{R}^{n+1}:\,\,\, \max ( |x|, |y|
)\le 1 \}
$$

For a natural $l ,\, 1\le l \le n+1$ let $\lambda_l (\xi , N)$ be
the $l$-th successive minimum of   ${\cal U}$ with respect to  $
\Lambda (\xi, N)$ and let $\mu_l (\xi , N)$ be the $l$-th
successive minimum of   ${\cal W}$ with respect to   $ \Lambda
(\xi, N)$. By the second Minkowski   theorem  for convex body (see
\cite{Cassil}, Ch. VIII or \cite{SCH}, Ch.IV)  we have
$$
\frac{1}{(n+1)!}  \le \prod_{l=1}^{n+1} \lambda_l (\xi , N)
 \le 1 ,\,\,\,
\frac{2^{n}}{w_n(n+1)!}  \le \prod_{l=1}^{n+1} \mu_l (\xi , N)
   \le \frac{2^{ n}}{w_n}.
$$

In the case $n=1$ we have $\mu_l (\xi , N)=\lambda _l (\xi , N), l
= 1,2$ .
 Suppose that  $\xi_1/\xi_2\not\in \mathbb{Q}$. Since  there are arbitrary large values of $N$ with
$\mu_1(\xi , N)=\mu_2 (\xi , N)$,   it may never happen that $
\mu_1(\xi , N)\to 0,\, N\to +\infty$.

Consider the general case. Suppose that  numbers $\xi_1,...,\xi_n$
are linearly independent over $\mathbb{Z}$ together with $1$. Then
for  $1\le k \le n$ there exist arbitrary large values of $N$ such
that  $\mu_k(\xi , N)=\mu_{k+1} (\xi , N)$. But in the case $n>1$
it may happen that $ \mu_1(\xi , N)\to 0,\, N\to +\infty$.
Moreover from A.Khinthcine's  result \cite{HINS} (see also
\cite{DS}) it follows that it  may happen that  $ \mu_{n-1}(\xi ,
N)\to 0,\, N\to +\infty$.

In this paper we prove the following theorem.

{\bf Theorem 1.}\,\,\,{\it Let $ 1\le k \le n-1$. Then there exist
real numbers $\xi_j \in[0,1), 1\le j\le n$, such that

$\bullet$\,\,   $ 1,\xi_1,...,\xi_n$ are linearly independent over
$\mathbb{Z}$;

$\bullet$\,\, $\mu_k (\xi, N) \to 0$ as $ N\to \infty$;

$\bullet$\,\, $\mu_{k+2} (\xi, N) \to \infty$ as $ N\to \infty$. }

We  make two remarks.

{\bf Remark 1.} The analogous  result for  $\lambda_l (\xi , N)$
was conjectured by W.M. Schmidt in \cite{S}. In this paper we
consider the Euclidean norm only for simplicity reasons. We must
note that the main result is valid not only for the Euclidean norm
but also for the sup-norm $|\cdot |_s$ (as it was conjectured in
W.M. Schmidt's paper \cite{S}).

{\bf Remark 2.} It is shown in Section 3 that Theorem 1 becomes
trivial without the condition on  $ 1,\xi_1,...,\xi_n$ to be
linearly independent over $\mathbb{Z}$.

The construction in the proof of Theorem 1 in the case $k=1$ is
very simple. It
 is close to the construction from   \cite{MUMN}, where the
author gives a counterexample to J. Lagarias' conjecture
concerning the behavior of  consecutive best simultaneous
Diophantine approximations (see \cite{Lag2}). We give a complete
proof of Theorem 1 in the case $ k=1$ in Section 2.

In the case $k>1$ the construction in the proof of Theorem 1  is a
little bit more difficult.  It is  close to procedures from
\cite{AM},\cite{MMZ} (See the author's review \cite{ME} for
related topics). We give a complete proof of Theorem 1 in the case
$ k >1$ in Sections 3-5.

In the proof we shall need the following notation.

 By
$\mu_l ({\cal C};L )$ we denote the $l$-th successive minimum of a
convex 0-symmetric  set ${\cal C}$ with respect to a lattice $L$.

Let $w_l$ denote the volume of the unit ball in the
$l$-dimensional Euclidean space.

For a set ${\cal M} \subset \mathbb{R}^{n+1}$     we denote by
$\overline{\cal M} $ the smallest closed set containing ${\cal
M}$. We also denote the smallest linear and affine subspaces of
$\mathbb{R}^{n+1}$ containing
 ${\cal M}$, by
 ${\rm span}\,{\cal M} $ and ${\rm
aff}\,{\cal M} $, respectively.

Consider a sublattice $\Lambda \subset \mathbb{Z}^{n+1}$.  By
${\rm dim}\,\Lambda$ we denote the dimension of the linear
subspace ${\rm span}\, \Lambda$. A sublattice $\Lambda\subset
\mathbb{Z}^{n+1}$ is defined to be {\it complete} if
$$
\Lambda = ({\rm span}\,\Lambda )\cap \mathbb{Z}^{n+1},
$$
that is  in the linear subspace ${\rm span}\, \Lambda$ there is no
integer points  different from points of $\Lambda$.

For every positive $Q$ and   $\sigma$ we define a cylinder ${\cal
C}_\xi (Q, \sigma) \subset \mathbb{R}^{n+1}$ as follows:
$$
{\cal C}_\xi (Q, \sigma) = \left\{ z=(x,y) \,\, x\in
\mathbb{R},\,\, y =(y_1,...,y_n) \in \mathbb{R}^n :\,\,\, |x|<Q,
\left( \sum_{j=1}^n|y_j-\xi_j x|^2\right)^{1/2}<\sigma \right\}.
$$
The quantities $ \mu_l(\xi , N)$ coincide with  the successive
minima of  $\mathbb{Z}^{N+1}$  with respect to $ {\cal C}_\xi (N,
N^{-1/n} )$, that is
$$
\mu_l(\xi , N) = \mu_l ({\cal C}_\xi (N, N^{-1/n}
);\mathbb{Z}^{n+1}).
$$

\section{ Proof of Theorem 1: case $ k=1$}
 We need  two auxiliary
results - Lemmas A and B.

{\bf Lemma A.}\,\,\,{\it Let $\xi = \left( 1, \frac{a_1}{q}, ... ,
\frac{a_n}{q}\right)\in [0,1]^{n+1}$ be a rational vector. Suppose
that for integers $q,a_1,...,a_n \in \mathbb{Z}$ we have
$$ q\ge 1, (q,a_1,...,a_n) = 1.
$$
Then for any positive $U >0$ and   any natural $i$ there exists a
positive  real number
$$
\eta = \eta (\xi, i, U) >0
$$
such that for every real vector $\xi' = (1, \xi_1',...,\xi_n')$
under condition $ |\xi' - \xi| < \eta $ the inequalities
$$
\mu_1 (\xi', N) \le i^{-1},\,\,\,\, \mu_2(\xi', N) \ge i$$ are
valid for all  $N$ in the interval
$$
(2qi\sqrt{n+1})^n \le N \le U.
$$}

Proof. First of all we note that  for $N \ge q$ we obviously have
$$
\mu_1 (\xi, N) \le qN^{-1}.
$$
Besides that, the Euclidean distance between the  one-dimensional
subspace $ {\rm span} \, \xi $
 and the set $ \mathbb{Z}^{n+1}\setminus {\rm span} \, \xi $ is
 not less than $ (q\sqrt{n+1})^{-1}$. Thus, in order to catch  an integer point,
independent with $\xi$,
  in the
 cylinder $t{\cal C}_\xi (N, N^{-1/n})$, we should take $t$ to be not less than $ N^{1/n}(q\sqrt{n+1})^{-1}$. Hence
 $$ \mu_2(\xi , N ) \ge N^{\frac{1}{n}}q^{-1}.
 $$
 From the hypothesis of Lemma A we deduce that the inequalities
 $$
\mu_1 (\xi, N) \le (2i)^{-1},\,\,\,\, \mu_2(\xi, N) \ge 2i$$ hold
for all $N \ge (2qi\sqrt{n+1})^n$. Now Lemma A follows from the
observation that
 for any $l$ the function $ \mu_l (\xi , N)$ is a continuous function in $\xi$ and $N$.
 Lemma A is proved.

{\bf Lemma B.}\,\,\,{\it Let $\Gamma $ be a two-dimensional
complete sublattice of $\mathbb{Z}^{n+1}$. Let $R$  be the
two-dimensional fundamental volume of $\Gamma$ and let $\rho
=\rho(\Gamma ) >0$ be the Euclidean distance between   $ {\rm span
} \,\Gamma$ and $\mathbb{Z}^{n+1}\setminus \Gamma$.
 Suppose
that $ \xi = (1,\xi_1,...,\xi_n)\in {\rm span } \,\Gamma$. Then
for any positive $N$ we have the following estimates:
$$
\mu_1(\xi ,N ) \le N^{\frac{1-n}{2n}} R^{\frac{1}{2}},\,\,\, \mu_3
(\xi, N) \ge N^{\frac{1}{n}} \rho .
$$}

Proof.  First of all we prove the upper bound for $\mu_1(\xi ,N
)$.
 The intersection of the cylinder
${\cal C}_\xi (N, N^{-1/n})$ with   $ {\rm span } \,\Gamma$ is an
$0$-symmetric parallelogram,whose two-dimensional volume greater
than or equal to $4N^{\frac{n-1}{n}}$. Suppose that $
4t^2N^{\frac{n-1}{n}} > 4R$ for some $t >0$. Then, by the
Minkowski convex body theorem   there is a nonzero point of
$\Gamma$  inside the parallelogram $t{\cal C}_\xi (N,
N^{-1/n})\cap {\rm span } \,\Gamma$. So, for any  $ t
> N^{\frac{1-n}{2n}} R^{\frac{1}{2}}$ the cylinder $t{\cal C}_\xi (N, N^{-1/n})$ contains a nonzero integer point   and the upper bound for
$\mu_1(\xi ,N )$ is proved.

To prove the lower bound for
 $\mu_3
(\xi, N) $ we need to take into account that  if
 the cylinder $t{\cal C}_\xi
(N, N^{-1/n})$ contains more than two linearly independent integer
points, then one of these points does not belong to   $ \Gamma$
and
 $ tN^{-1/n}\ge \rho$.
Lemma B is proved.

Now we describe the inductive procedure which  gives the proof of
Theorem 1 in the case $k=1$.

The set of all $n$-dimensional sublattices of $\mathbb{Z}^{n+1}$
is countable. We fix an enumeration of this set and  let
$$ L_1, L_2,..., L_i , ...
$$
be the sequence of all $n$-dimensional complete sublattices of
$\mathbb{Z}^{n+1}$.
 Set $\pi_i = {\rm span} \,L_i$. Suppose that $$\pi_1 = \{z=(x,y_1,...,y_n)\in\mathbb{R}^{n+1}:\,\ x=0\}.$$

We construct a sequence of rational vectors
$$
\xi_i = \left( 1, \frac{a_{1,i}}{q_i}, ... ,
\frac{a_{n,i}}{q_i}\right),\,\,\,\, q_i,a_{1,i},...,a_{n,i} \in
\mathbb{Z},\,\,\, q_i\ge 1, (q_i,a_{1,i},...,a_{n,i}) = 1,\,\,\, i
= 1,2,3,... $$ with $ q_i \to \infty$ as $ i \to \infty$,
 a sequence of
two-dimensional complete sublattices $\Gamma_{i+1} \subset
\mathbb{Z}^{n+1},\,\, i=1,2,3,...$,
  and a sequence of positive real numbers
$$
\eta_1,\eta_2,...,\eta_i, ...
$$
satisfying the following conditions (i) - (iv).

(i) For every $ i \ge 1$ we have
$$
\xi_i, \xi_{i+1} \in {\rm span}\, \Gamma_{i+1}.
$$

(ii) The closed ball  $\overline{\cal B}_i$ of radius $\eta_i$
centered at $\xi_i$ and   has no common points with the subspace
$\pi_i$:
$$
\overline{\cal B}_i\bigcap \pi_i = \varnothing.
$$

(iii)  The   balls $\overline{\cal B}_i$   form a nested sequence:
$$
\overline{\cal B}_1 \supset \overline{\cal B}_2 \supset ...
\supset \overline{\cal B}_i .
$$

(iv) Let $H_0=1$ and
$$
 H_i = (4q_i (i+1)\sqrt{n+1})^n,\,\, i = 1,2,3,... \,\, .
$$
Then for any   $ i \ge 1$, for any $\xi \in\overline{\cal B}_i$
and for any $N$, such that
$$ H_{i-1}\le N <H_i,
$$
the following inequalities holds:
$$
\mu_1(\xi, N) \le i^{-1},\,\,\, \mu_3 (\xi, N) \ge i.
$$

Suppose all the objects are already  constructed. Then from (ii),
(iii) one can easily see that $\lim_{i\to \infty} \eta_i = 0$.
Then the intersection $\cap_{i\in \mathbb{N}} \overline{\cal B}_i$
consists of the only one point. Note that for every $i$ the center
$\xi_i$ of the ball $\overline{\cal B}_i$  has its first
coordinate equal to one. So the unique point from the intersection
$\cap_{i\in \mathbb{N}} \overline{\cal B}_i$ is of the form
$$
  \xi
=(1,\xi_1,...,\xi_n)   $$ (as $\xi_i \to \xi,\,\, i \to \infty$).
Then
$$
\mu_1(\xi, N) \le i^{-1},\,\,\, \mu_3 (\xi, N) \ge i,\,\,\,
H_{i-1}\le N< H_i.
$$
Hence
$$
\mu_1(\xi, N) \to 0,\,\,\, \mu_3 (\xi, N) \to \infty,\,\,\, N\to
\infty,
$$
and it follows from the conditions (ii) and (iii)   that
$1,\xi_1,...,\xi_n$ are linearly independent over $\mathbb{Z}$.
This proves Theorem 1 in the case $k=1$.

We start our inductive procedure with the vector
$$
\xi_1 = (1, \underbrace {0,...,0}_{n\,\,\text{times}}) .$$ Then $
H_1 = (8\sqrt{n+1})^n$. The sublattice $\Gamma_1$ is not defined
yet and we do not   care about the condition (i) at this stage.
The condition (ii) is obviously  satisfied for any choice of $
\eta_1 < 1$. The conditions (iii) is   empty. Recall that $H_0 =
1$. To satisfy the condition (iv)   we take $\eta_1$ to be small
enough.

Now we pass to   the inductive step. Suppose that all the objects
$\xi_i, \Gamma_i, \eta_i $, $i \le t$ satisfying  conditions  (i)
-- (iv) are already constructed. We describe how to construct the
$(t+1)$-th set of objects.

First of all we can take a two-dimensional complete sublattice
$\Gamma_{t+1}$  satisfying the conditions
$$
q_t\xi_t = (q_t, a_{1,t},...,a_{n,t} ) \in \Gamma_{t+1},\,\,\,
\Gamma_{t+1} \not\subset \pi_{t+1} .
$$
Let $R_{t}$ be the two-dimensional fundamental volume of
$\Gamma_{t+1}$ and let $\rho_{t}$ be the Euclidean distance
between   $
 {\rm span } \,\Gamma_{t+1}$ and
$\mathbb{Z}^{n+1}\setminus \Gamma_{t+1}$.

Set
$$
U_{t} = \max\left( (2(t+1)\rho_{t}^{-1})^n,\,
(2(t+1)R_t^{1/2})^{\frac{2n}{n-1}} \right).
$$
Now we apply Lemma A with $\xi = \xi_t, i = 2(t+1)$ and $ U =
U_t$. We get a positive $\eta_t'$, such that for every $\xi '$
under the condition $|\xi ' - \xi_t|< \eta_t'$ one has
$$
\mu_1 (\xi', N) \le (2t+2)^{-1},\,\,\,\, \mu_2(\xi', N) \ge 2t+2
$$
for every $N$ in the interval
$$
H_t = (4q_t(t+1)\sqrt{n+1})^n \le N \le U_t. $$

 Obviously, there is an integer point
$$
(q_{t+1},a_{1,t+1},...,a_{n,t+1}) \in \Gamma_{t+1}\setminus
\pi_{t+1}, \,\,\, q_{t+1} \ge q_t,\,\,\,
(q_{t+1},a_{1,t+1},...,a_{n,t+1}) = 1,
$$ such that for
$$
\xi_{t+1} = \left( 1, \frac{a_{1,t+1}}{q_{t+1}}, ... ,
\frac{a_{n,t+1}}{q_{t+1}}\right)
$$
we have
$$
|\xi_{t+1} - \xi_t | < \frac{\min ( \eta_t, \eta_t ')}{2}.
$$
Since $\xi_{t+1} \in \Gamma_{t+1}$, we can apply Lemma B with $\xi
= \xi_{t+1}, \Gamma =\Gamma_{t+1}$. This gives that for any $N$
under the condition
\begin{equation}
N \ge U_t \label{mai}
\end{equation}
one has
\begin{equation}
\mu_1 (\xi_{t+1}, N) \le (2(t+1))^{-1},\,\,\,\, \mu_3(\xi_{t+1},
N) \ge 2(t+1). \label{maine}
\end{equation}
But $|\xi_{t+1} - \xi_t|<\eta_t'$ and $\mu_3(\xi_{t+1}, N)
\ge\mu_2(\xi_{t+1}, N) $. So, the inequalities (\ref{maine}) are
valid not only for $N$ in the interval (\ref{mai}) but also for
$N$ in the interval $ N \ge H_t$. Having constructed   $
\xi_{t+1}$, we define  $H_{t+1} $ from the condition (iv) of the
$(t+1)$-th step of the inductive process.
 Now we
take into account that
 for any $l$ the function $ \mu_l (\xi , N)$ is a continuous function in $\xi$ and $N$.
This means that we can find a number $\eta_{t+1} < {\min ( \eta_t,
\eta_t ')}/{2}, $ such that
$$
\mu_1 (\xi , N) \le (t+1)^{-1},\,\,\,\, \mu_2(\xi, N) \ge t+1
$$
for all $\xi$ under the condition
$$
|\xi - \xi_{t+1}| \le  \eta_{t+1}
$$
and  all $N$ in the interval
$$
H_t \le N < H_{t+1}.
$$
Moreover, since $\xi_{t+1}\not\in \pi_{t+1}$, we can take
$\eta_{t+1}$ to be small enough, so that the ball $\overline{\cal
B}_{t+1}$ of  radius $ \eta_{t+1}$ centered at $ \xi_{t+1}$ and
has no common points with  $\pi_{t+1}$. The $(t+1)$-th step of the
inductive procedure is described completely and hence Theorem 1 in
the case $ k = 1$ is proved.

\section{Lemmas concerning   successive minima and badly
approximable numbers}

In this section we start the proof of Theorem  1 in the general
case.

Let $l$ be an integer and  $ 1\le l \le n$. First of all we should
say that everywhere in the sequel we consider $l$-dimensional
lattices $L$ such that for some real vector $\xi =
(1,\xi_1,...,\xi_n) \in \mathbb{R}^{n+1}$ one has
$$
\xi \in {\rm span }\, L.
$$
This condition leads to the following corollary concerning linear
subspace ${\rm span }\, L $. Consider affine subspace
$${\cal P} =\{
z=(x,y_1,...,y_n)\in \mathbb{R}^{n+1}:\,\, x = 1\} \subset
\mathbb{R}^{n+1}.
$$
Then the intersection
$$
{\rm span }\, L\cap {\cal P}
$$
is an affine subspace of dimension $ l-1$.

 We
prove some auxiliary results on successive minima and badly
approximable numbers.

For positive integer $l$ and positive $\sigma \in (0,1)$ define
 \begin{equation}\label{starr}
Q_1= Q_1(\sigma, l,s) = 2^{l-1}w_{l-1}^{-1} \sigma^{1-l}s,\,\,\,
Q_2 =Q_2(\sigma, l,s) = 2^{l-1}Q_1 = 2^{2(l-1)}w_{l-1}^{-1}
\sigma^{1-l}s.
\end{equation}

{\bf Lemma 1.}\,\,\,{\it  Consider an integer  $l$, such that $
2\le l\le n+1$. Consider a complete sublattice $ L \subseteq
\mathbb{Z}^{n+1}$ and suppose that $ {\rm dim}\,  L  = l$  and $
\xi \in {\rm span}\, L$. Suppose that the fundamental
$l$-dimensional volume of the lattice $L$ is equal to $s$.
Consider a cylinder ${\cal C} = {\cal C}_\xi (Q, \sigma)$. Suppose
also that for some positive $Q$ and $\sigma $ we have $ {\cal C}
\bigcap L = \{ 0\}. $ Then the following upper bounds are valid:
$$
\mu_1 ({\cal C} )\le Q_1Q^{-1}, \,\,\,\, \mu_m ({\cal C} )\le  Q_2
Q^{-1} ,\,\,\, 2\le m \le l.
$$}

Proof. Consider the cylinder
$$
{\cal C}^{(1)}= {\cal C}_\xi (Q_1,\sigma)\bigcap {\rm span }\,L.
$$
As $ \xi \in {\rm span }\,L$ we see that each section
$$
{\cal C}^{(1)}\cap \{z=(x,y_1,...,y_n) \in \mathbb{R}^{n+1}:\,\, x
= x_0\}, \,\,\, |x_0| \le Q_1
$$
is a $(l-1)$-dimensional ball of the volume $w_{l-1}\sigma^{l-1}$.
Let $ H_1$ be the distance between $(l-1)$-dimensional facets of
the cylinder ${\cal C}^{(1)}$.
 Then
 the $l$-volume
of $ {\cal C}^{(1)}$ is equal to
$$
w_{l-1} H_1 \sigma^{l-1}.$$ But $H_1 \ge 2Q_1$. So the $l$-volume
of $ {\cal C}^{(1)}$ is
$$
\ge 2w_{l-1} Q_1 \sigma^{l-1} = 2^ls.
$$
(in the equality here we take into account the definition of $Q_1$
from (\ref{starr})).
 By the
Minkowski's theorem for  convex body, there is a nonzero  integer
points $\zeta^{(1)} \in {\cal C}^{(1)}\cap L$.
 So, $\mu_1 ({\cal C} )\le Q_1Q^{-1}$
and the bound for the first successive minimum is proved.

Here we should note that as $\sigma <1$ the first coordinate of
$\zeta^{(1)}$ is not equal to zero.

 Now we  describe an inductive  process of constructing
linearly independent integer point $ \zeta^{(1)} ,...,
\zeta^{(l)}$ with non-zero first coordinates which ensure the
upper bound for the successive minima under consideration.

 Suppose that linearly independent integer points $ \zeta^{(1)} ,..., \zeta^{(\nu)}\in  {\rm span }\, L$ with $1\le \nu
\le l-1$ are already constructed. Set $ \pi = {\rm span }
(\zeta^{(1)} ,..., \zeta^{(\nu)})\subset {\rm span } \, L$.   As
all the points $ \zeta^{(1)} ,..., \zeta^{(\nu)}$ are linearly
independent we see that ${\rm dim }\,\pi = \nu < l$. Note that the
dimension of the affine subspace $  \pi' = \pi\cap
\{z=(x,y_1,...,y_n) \in \mathbb{R}^{n+1}:\,\, x = Q\}$ is equal to
$ {\rm dim }\,\pi-1< l-1  $.
  Consider the  facet $ {\cal B} = \overline{\cal
C}\cap \{ x = Q\}$. This facet is an $n$-dimensional ball of
radius $\sigma$ centered at $Q \xi$.

Note that as $ \xi \in {\rm span}\, L$ we see that the
intersection ${\rm span}\,L\cap\{ x = Q\}$ is a
$(l-1)$-dimensional affine subspace. Consider the intersection $
{\cal B}\cap {\rm span}\, L$. As $Q \xi\in {\rm span}\, L $ we see
that this intersection is a $(l-1)$-dimensional ball centered at
$Q \xi$.

 We have the following situation.  In the $(l-1)$-dimensional affine
subspace  ${\rm span}\,L\cap\{ x = Q\}$ there are the ball $ {\cal
B}\cap {\rm span}\, L$ of dimension $l-1$ and the affine subspace
$ \pi'  $ of dimension ${\rm dim }\, \pi'   < l-1$. So there
exists a $n$-dimensional ball ${\cal B}'\subset {\cal B} \subset
\{ x = Q\}$ of radius $\sigma/2$ centered at a certain point
 $\Xi \in {\rm span}\,L\cap\{ x = Q\}$  and such  its $(l-1)$-dimensional section $ {\cal B}'\cap {\rm span}\, L$ does not intersect with
 $\pi'$:
$$ {\cal B}'\cap {\rm span}\,L \cap \pi =\varnothing.$$
In fact as $ \pi' \subset {\rm span}\,L \cap\{ x = Q\}$, it means
that
$$ {\cal B}'\cap \pi =\varnothing.$$
Put $ \xi^{(\nu+1)} = \frac{\Xi}{Q}$ and consider the cylinder
$$
{\cal C}^{(\nu+1)}= {\cal C}_{\xi^{(\nu+1)}} (Q_2,\sigma /2
)\bigcap {\rm span}\, L .
$$
As $\xi^{(\nu+1)}\in  {\rm span}\, L $ from (\ref{starr}) we see
that the
  $l$-volume
of $ {\cal C}^{(\nu+1)}$ is equal to
$$
w_{l-1} H_2\left(\frac{\sigma}{2}\right)^{l-1}  ,$$ where  $
H_2\ge 2Q_2$ is the distance between $(l-1)$-dimensional facets of
the cylinder ${\cal C}^{(\nu+1)}$. So the
  $l$-volume
of $ {\cal C}^{(\nu+1)}$ is
$$
\ge 2w_{l-1} Q_2\left(\frac{\sigma}{2}\right)^{l-1} =
 2^{l}s.
 $$
  Applying again  the Minkowski theorem for convex body   we get a nonzero  integer point $\zeta^{(\nu+1)} \in
{\cal C}^{(\nu+1)}\cap L$.

As $\sigma <1$ we see that the first coordinate of
$\zeta^{(\nu+1)}$ is not equal to zero. Moreover the first
coordinate of the point $\zeta^{(\nu+1)}$ is greater than $Q$ as
there is no nonzero integer points in ${\cal C} \bigcap {\rm
span}\, L$ and
$$
{\cal C}^{(\nu+1)} \cap \{z=(x,y_1,...,y_n) \in
\mathbb{R}^{n+1}:\,\, |x| \le  Q\} \subset  {\cal C} =  {\cal
C}_\xi (Q, \sigma).$$ So
$$
\zeta^{(\nu+1)}\not \in {\cal C}.
$$
But
$$
\pi \cap {\cal C}^{(\nu+1)} \subset  {\cal C}.
$$
So $\zeta^{(\nu+1)}\not \in \pi$. It means that $\zeta^{(\nu+1)}$
is independent of  $ \zeta^{(1)} ,..., \zeta^{(\nu)}$.

To conclude the proof we make two following observations:

1. each point $z = (x,y_1,...,y_n)\in {\cal C}^{(\nu+1)}$
satisfies $ |x| \le Q_2$;

2. for the section $ {\cal C}^{(\nu+1)}\cap \{ x = Q\}$ one has
$$
{\cal C}^{(\nu+1)}\cap \{ x = Q\} \subset {\cal C} \cap \{ x = Q\}
,
$$
and the section ${\cal C} \cap \{ x = Q\}$ is a ball of radius
$\sigma$ centered at $Q\xi$.

 So
$$
{\cal C}^{(\nu+1)}\,\, \subset \,\, {\cal C}_{\xi^{(\nu+1)}}
(Q_2,Q_2Q^{-1}\sigma).$$ Now $ \mu_m ({\cal C} )\le Q_2Q^{-1} $
for $ 2\le m \le l. $ Lemma 1 is proved.

{\bf Lemma 2.}\,\,\,{\it  Let $ 2\le l\le n+1$. Consider a
complete sublattice $ L \subseteq \mathbb{Z}^{n+1}$ and suppose
that $ {\rm dim}\,  L  = l$  and $ \xi \in {\rm span}\, L$.
Suppose that the fundamental $l$-dimensional volume of  $L$ is
equal to $s$. Suppose also that for some $Q ,\sigma >0 $ we have
$$
{\cal C}_\xi (Q, \sigma)\bigcap L = \{ 0\}.
$$
Then for any  $M , \delta >0$ the following upper bound is valid:
$$
\mu_{l}({\cal C}_\xi (M,\delta) )\le  Q_2 Q^{-1}  \max( QM^{-1},
\sigma \delta^{-1}).
$$
 }

{\bf Corollary.}\,\,\,{\it  Suppose that the conditions of Lemma 2
are satisfied. Then for the cylinder ${\cal C}_\xi (N,N^{ -1/n})$
we have
$$
\mu_{l}(\xi , N)\le  Q_2 Q^{-1} \max( QN^{-1}, \sigma N^{ 1/n}).
$$}

Proof of Lemma 2. Put $ t = \max( QM^{-1}, \sigma \delta^{-1})$.
Then
$$
{\cal C}_\xi (Q, \sigma)\subset t {\cal C}_\xi (M,\delta),
$$
and applying Lemma 1 we see that
$$
\mu_{l}({\cal C}_\xi (M,\delta) ) = t \mu_{l }(t{\cal C}_\xi
(M,\delta) )\le t\mu_{l}({\cal C}_\xi (Q, \sigma))\le  Q_2 Q^{-1}
t.
$$
Lemma 2 is proved.

Put $ \xi_0 = 1$. For a real vector $\xi = (\xi_0,\xi_1,...,\xi_n
)=(1,\xi_1,...,\xi_n ) \in\mathbb{R}^{n+1}$ we define ${\rm
dim}_\mathbb{Q}\xi $ to be the maximal integer  $t$, such that the
components $ \xi_{j_1},...,\xi_{j_t}, 0\le j_1,...,j_t \le n+1$
are linearly independent over $\mathbb{Q}$. For example, the
equality ${\rm dim}_\mathbb{Q}\xi = 1$ occurs only if $\xi \in
\mathbb{Q}^{n+1}\setminus \{ 0\}$ and the equality ${\rm
dim}_\mathbb{Q}\xi = n+1$ occurs only if all the components
$1,\xi_1,...,\xi_n$ are linearly independent over $\mathbb{Q}$.
Obviously, if ${\rm dim}_\mathbb{Q}\xi = l, \,\,\, 1\le l \le
n+1$, then there is a complete  sublattice $ L \subseteq
\mathbb{Z}^{n+1}$, such that
  $
  {\rm dim}\,  L = l  $ and $ \xi \in {\rm span}\, L$. Moreover,
$$
{\rm dim}_\mathbb{Q}\xi =\min\{ l\in \mathbb{N}:\,\, \text{there
exists a sublattice} \, L \subseteq \mathbb{Z}^{n+1}, \,\,\text{
such that}\,\,
  {\rm dim}\,  L  = l\,\,\text{ and}\,\,
 \xi \in {\rm span}\, L \}.$$

Let us now we consider a complete sublattice $ L \subseteq
\mathbb{Z}^{n+1}$,  such that
  $
  {\rm dim}\, L   = l\ge 2 $  and let us consider a vector
  $ \xi =(1,\xi_1,...,\xi_n )\in {\rm span}\, L$ (then $ {\rm
dim}_\mathbb{Q}\xi\le l$). We shall say that
 $\xi $ is   {\it $\gamma$-badly approximable with
respect to  } $L$ (briefly $(L,\gamma)$-BAD) if for any nonzero
integer point $ \zeta =(q,a) = (q,a_1, ...,a_n) \in L $ with
$q\neq 0$ one has
\begin{equation}
|q\xi - \zeta| \ge \gamma |q|^{-1/(l-1)}. \label{BAD}
\end{equation} We should note that for any $(L,\gamma)$-BAD vector
$\xi$ and any $Q\ge 1$ the cylinder
\begin{equation}
{\cal C}_\xi (Q,  \sigma_Q)\bigcap {\rm span}\, L,\,\,\, \sigma_Q
= \gamma Q^{-1/(l-1)} \label{CY} \end{equation} contains no
nonzero integer points inside. A vector $\xi \in {\rm span}\, L$
is defined to be {\it badly approximable with respect to  } $L$
(briefly $ L$-BAD) if (\ref{BAD}) holds with some positive
$\gamma$. It is easy to see that if vector $\xi$ is badly
approximable with respect to   $L$ and ${\rm dim}\,  L   = l$ then
$\xi \in {\rm span}\,\Lambda$ and  $ {\rm dim}_\mathbb{Q}\xi = l$.

Let $W\ge 1$. It is necessary for us  to consider vectors $\xi \in
{\rm span} \, L$, such that  the cylinder (\ref{CY})  contains no
nonzero integer point inside only for $ Q\ge W$. We define such
vectors  to be {\it $(\gamma, W)$-badly approximable with respect
to } $L$ (briefly $(L,\gamma, W)$-BAD). A vector $\xi \in {\rm
span} L$ is  $(\gamma, W)$-badly approximable with respect to
$L$ iff (\ref{BAD}) holds for all $\zeta $ with $|q|\ge W$ and for
all $q$ under the condition $1\le |q|\le W$   the following
inequality holds instead of (\ref{BAD}):
$$
|q\xi - \zeta| \ge \gamma W^{-1/(l-1)}.
$$
It is obvious that a vector $\xi \in {\rm span} \, L$ is
$(\Lambda,\gamma)$-BAD iff it is $(\Lambda,\gamma, 1)$-BAD.

{\bf Example 1.}\,\,\, Consider the space $\mathbb{R}^{n+1} $
related to   coordinates $x,y_1,...,y_n$. Consider the case when
real algebraic integers $1,\alpha_1,...,\alpha_{l-1}$ form a basis
of a real algebraic field ${\cal K}$ of degree $l\ge 2$. Then
there exists a constant $\gamma = \gamma ({\cal K}) $, such that
for all natural $q$ we have
$$
\left( \sum_{j=1}^{l-1}||q\alpha_j||^2\right)^{1/2} \ge \gamma
q^{-1/(l-1)}
$$
(see \cite{Cas}, Chapter V, \S 3) and hence the
$(n+1)$-dimensional vector
$$
(1,\alpha_1,...,\alpha_{l-1},\underbrace
{0,...,0}_{n+1-l\,\,\text{times}})
$$
is $(L, \gamma ({ \cal K}))$-BAD  where $L = \mathbb{Z}^{n+1}\cap
\{y_l = ...= y_n = 0\}$.

Define
$$
G_1= G_1(l,s,\gamma)=2^{2(l-1)}w^{-1}_{l-1}s\gamma^{1-l},\,\,\,
G_2= G_2(l,s,\gamma)=
2^{2(l-1)}w^{-1}_{l-1}s\gamma^{-\frac{(l-1)^2}{l}}.
$$

{\bf Lemma 3.}\,\,\,{\it   Let $ L \subseteq \mathbb{Z}^{n+1}$ be
a  complete sublattice, such that
  $
  {\rm dim}\,  L   = l \ge 2 $  and let
  $ \xi =(1,\xi_1,...,\xi_n )\in {\rm span} \, L$ be an $(L,\gamma , W)$-BAD
  vector.  Let $s$ be the $l$-dimensional fundamental volume of   $L $.
  Consider  positive $M,\delta $ and the cylinder
${\cal C}= {\cal C}_\xi (M,\delta)$. Then the following statements
hold.

1) If
\begin{equation}
  \left(M\gamma \delta^{-1}\right)^{\frac{l-1}{l}}\le W
\label{EN0B}
\end{equation}
    then
  \begin{equation}
  \mu_l ({\cal C}) \le G_1 WM^{-1}.
\label{C1C0B}
\end{equation}

  2) If
\begin{equation}
\left(M\gamma \delta^{-1}\right)^{\frac{l-1}{l}}\ge W \label{ENB}
\end{equation}
   then
  \begin{equation}
  \mu_l ( {\cal C}) \le G_2
M^{-\frac{1}{l}} \delta^{\frac{1-l}{l}}
  .
\label{C1CB}
\end{equation}
  }

 {\bf Remark.}\,\, We actually construct  in the proof
   $l$ nonzero linearly independent integer points $\zeta_j \in L$ lying in
 the  cylinder $\mu_l ({\cal C}) \cdot   \overline {\cal C}$.
 It is seen from the construction that in the
 case 2) of Lemma 3 each ray $[0,\zeta_j)$, $ 1\le j \le l$
 intersects the facet  $\{ x = M\}$ of the cylinder ${\cal C}= {\cal C}_\xi (M,\delta)$.

  Proof of Lemma 3.
For $Q\ge W$ the cylinder (\ref{CY}) has no nonzero integer
points. By  Lemma 2 for any $ Q\ge W$ we have
$$
\mu_{l}( {\cal C})\le G_1 \max( QM^{-1}, \gamma Q^{-1/(l-1)}
\delta^{-1}).
$$
Consider
$$
m(M,\delta , W) =\min_{Q\ge W} \max( QM^{-1}, \gamma Q^{-1/(l-1)}
\delta^{-1}) .
$$
If (\ref{EN0B}) holds we have
$$
m(M,\delta , W)= WM^{-1}.
$$
If (\ref{ENB}) holds we see that
$$
m(M ,\delta, W)=  \gamma^{\frac{l-1}{l}}\delta^{\frac{1-l}{l}}
M^{-\frac{1}{l}} .
$$
Lemma 3 follows.

Lemma 3 applied to the cylinder ${\cal C}_\xi (N, N^{-1/n})$ gives
the following

 {\bf Corollary 1.}\,\,\,{\it Let $\xi =(1,\xi_1,...,\xi_n )
\in \mathbb{R}^{n+1}$. Let $ L \subseteq \mathbb{Z}^{n+1}$ be a
complete sublattice such that
  $
  {\rm dim}\,  L = l \ge 2 $  and let
  $ \xi \in {\rm span}\, L$ be a $(L,\gamma , W)$-BAD
  vector.  Let $s$ be the $l$-dimensional fundamental volume of   $L $.
  Then

1) for any positive $N$ under the condition
\begin{equation}
   N \le \gamma^{-\frac{n}{n+1}} W^{\frac{ln}{(n+1)(l-1)}}
\label{EN0}
\end{equation}
    one
  has
  \begin{equation}
  \mu_l (\xi , N) \le G_1   WN^{-1}.
\label{C1C0}
\end{equation}

  2) for any  $N$ under the condition
\begin{equation}
   N \ge \gamma^{-\frac{n}{n+1}} W^{\frac{ln}{(n+1)(l-1)}}
\label{EN}
\end{equation}
    one
  has
  \begin{equation}
  \mu_l (\xi , N) \le G_2N^{\frac{l-n-1}{nl}}.
\label{C1C}
\end{equation}
  }

{\bf Corollary 2.}\,\,\,{\it Let $ 2\le l \le n$. Let $\xi \in
\mathbb{R}^{n+1}$. Let $ L \subseteq \mathbb{Z}^{n+1}$ be a
complete sublattice  such that
  $
  {\rm dim}\, L   = l  $  and let
  $ \xi \in {\rm span}\, L$ be a $L$-BAD
  vector.Then
  $$
  \mu_l (\xi , N) \to 0,\,\,\,\mu_{l+1} (\xi , N) \to +\infty,\,\,\, N \to \infty .
  $$}

{\bf Remark 1.}\,\,\, Obviously, for $l = 1$ in the case $ {\rm
dim}_\mathbb{Q}\xi = 1$ we have $$
  \mu_1 (\xi , N) \to 0,\,\,\,\mu_{2} (\xi , N) \to +\infty,\,\,\, N \to \infty .
  $$

{\bf Remark 2.}\,\,\ Of course, the assertion of Corollary 2
enforces the components $1,\xi_1,...,\xi_n$ to be linearly
dependent over $\mathbb{Q}$.

  Proof of Corollary 2. The statement about $\mu_l (\xi , N)$ follows  immediately
  from Corollary 1 of
  Lemma 3
  as the exponent in the right hand side of (\ref{C1C}) is negative.  We
  prove the statement about $\mu_{l+1} (\xi , N)$.

Suppose that $f_1, ...,f_l\in L $ form a basis of   $L$. Then it
can be completed to a basis $f_1,...,f_l, g_{l+1},...,g_{n+1}$ of
the entire integer lattice $\mathbb{Z}^{n+1}$. Let $L '$ be the
sublattice generated by $g_{l+1},...,g_{n+1}$. Then $
\mathbb{Z}^{n+1}= L \oplus L '$, \,\,$ {\rm dim }\, L '   = n+1-l$
and $ {\rm span }\,L \cap {\rm span }\,L ' =\{ 0\}.$ Consider the
$n+1-l$ dimensional linear subspace $\pi \subset \mathbb{R}^{n+1}
$, orthogonal to ${\rm span }\,L $. Then $\pi \oplus {\rm span
}\,L= \mathbb{R}^{n+1}$ and any two vectors $ u\in \pi,v \in {\rm
span }\,L$ are orthogonal. Hence the orthogonal projection of $L'$
onto   $\pi$ is a lattice $L ''$, such that ${\rm span}\, L '' =
\pi$. Let $\omega=\omega (L )>0 $ be the length of the shortest
nonzero  vector in $L ''$. Then for any integer point $\zeta \in
\mathbb{Z}^{n+1} \setminus L$ the Euclidean distance from $\zeta $
to   $ {\rm span }\,L$ is not less than $   \omega$. Suppose that
$ \xi \in {\rm span }\,L$ and that the cylinder $ {\cal C}_\xi
(tN, tN^{-1/n} )$ contains  $ l+1$ linearly independent integer
points. Then at least one of these   points belongs to
$\mathbb{Z}^{n+1} \setminus L$. Hence $ tN^{-1/n}
> \omega (L )$ and $\mu_{l+1} (\xi , N) \ge \omega (L
)N^{1/n}\to +\infty,\,\, N\to \infty$. The Corollary  is proved.

{\bf Lemma 4.}\,\,\,{\it Let $ 2\le l \le n+1$. Let  $s$ be the $l
$-dimensional fundamental volume of a lattice $L$, ${\rm dim}\ L =
l .$ Suppose that a vector $\xi = (1,\xi_1,...,\xi_n) \in {\rm
span }\,L $ and positive numbers $\gamma  $ and $ T \ge 1$ satisfy
the equality
\begin{equation}
{\cal C}_{\xi } (T,\gamma  T^{-1/(l-1)}) \bigcap L = \{ 0\}.
\label{EMP}
\end{equation}
Let
\begin{equation} \gamma ^*= \gamma^*(\gamma,L) = \min \left( 3^{-2}\gamma ,
3^{-l-2} (2 w_{l-1}l !)^{-\frac{1}{l-1}} s^{\frac{1}{l-1}}
 \right) .
\label{prime}
\end{equation}

Then there exists an $(L,\gamma^*, T)$-BAD vector $\xi^*=
(1,\xi_1^*,...,\xi_n^*) \in {\rm span }\, L$, such that
\begin{equation}
|\xi^*-\xi  |< \gamma  T^{-\frac{l }{l-1}} \label{prim}.
\end{equation}  }

Proof. Put $T_\nu = 3^{ (l-1)\nu  }T,\,\, \nu = 0,1,2,... $. To
prove Lemma 4 it  suffices to construct a sequence  of cylinders
$$
{\cal C}^{(\nu )} = {\cal C}_{\xi^{(\nu ) }} (T_\nu,9\gamma^*
T_\nu^{-1/(l-1)}),\,\,\,\,\, \xi^{(\nu ) }\in {\rm span }\,L
$$
 such that

(i) for every $\nu$ we have ${\cal C}^{(\nu )} \cap L =\{ 0\}$;

(ii) the section $  \{ x = T_{\nu}\}\cap {\cal C}^{(\nu +1)} $ of
the cylinder ${\cal C}^{(\nu + 1)}$ lies inside the facet $ \{ x =
T_{\nu}\}\cap \overline {{\cal C}^{(\nu )}} $ of the preceding
cylinder ${\cal C}^{(\nu )} $; moreover, the distance between the
centers of   ${\cal B}$ and   ${\cal B}'$ does not exceed  $
3\gamma^* T_\nu^{-1/(l-1)} $.

If such cylinders ${\cal C}^{(\nu )} $ are constructed  and ${\cal
C}^{(0 )} = {\cal C}_{\xi } (T,\gamma T^{-1/(l-1)})$ then the
vector $ \xi ^* =\lim_{\nu \to + \infty} \xi^{(\nu )}$ satisfies
(\ref{prim}). Moreover,  it follows  from (ii) that $\xi ^*$ is an
$(L,\gamma^*, T)$-BAD vector. Indeed, for an integer point $\zeta
=(q,a_1,...,a_n)$ with $ T_\nu \le q\le T_{\nu+1} = 3^{l-1}T_\nu$
we have
$$
|q\xi^*-\zeta|\ge 3\gamma^*T_{\nu+1}^{-1/(l-1)}\ge
\gamma^*q^{-1/(l-1)}.
$$

 We now  describe the inductive process, which  constructs the sequence of
cylinders ${\cal C}^{(\nu )}$.  Suppose that ${\cal C}^{(\nu )}$
is already constructed. Consider the cylinder
$${\cal C}'=
{\cal C}_{\xi^{(\nu ) }} (T_{\nu+1},3^{l+1 }\gamma^*
T_\nu^{-1/(l-1)}).
$$
We prove that there exists a linear subspace ${\cal L} \subset
{\rm span} \, L $ of dimension ${\rm dim }\, {\cal L} = l-1$
containing  all the integer points $\zeta \in L \bigcap {\cal C} '
$. Suppose that there are $l $ linearly independent integer points
$$
\zeta^{(1)},..., \zeta^{(l)} \in   \Gamma\bigcap{\cal C} '  .
$$
Then the $l$-dimensional volume $V$ of  the convex hull ${\rm
conv}\,(0, \zeta^{(1)},..., \zeta^{(l)})$ is bounded from below by
the fundamental volume of $L$:
\begin{equation}
V \ge s ( l!)^{-1}. \label{contra}\end{equation} On another hand,
the volume of ${\rm conv}\,(0, \zeta^{(1)},..., \zeta^{(l)})$
admits an  upper bound  based on the relation  ${\rm conv}\,(0,
\zeta^{(1)},..., \zeta^{(l)} ) \subset {\cal C} ' . $ Taking into
account (\ref{prime}) we see that
\begin{equation}
V \le T_{\nu+1}\left( 3^{l+1}\gamma^*
T_\nu^{-1/(l-1)}\right)^{l-1} \le s (2 l !)^{-1}
  \label{contrb}\end{equation}
 Relations (\ref{contra},\ref{contrb})   contradict  each other, which means that all the  integer points from the cylinder under
 consideration lie in a subspace ${\cal L}$.

 Let ${\cal B}$ be the $(l-1)$-dimensional
facet $\{ x = T_\nu \}$ of ${\cal C}^{(\nu )}\cap {\rm span}\, L$.
In fact ${\cal B}$ is an $(l-1)$-dimensional open ball of  radius
$3\gamma^* T_\nu^{-1/(l-1)}$ centered at
   $T_\nu \xi^{(\nu )}\in   {\rm span}\, L$. There is an $(l-1)$-dimensional open ball ${\cal B}'\subset {\cal B}$ of radius
   $\gamma ^* T_\nu^{-1/(l-1)} = 3\gamma ^* T_{\nu+1}^{-1/(l-1)}$ and centered at a certain point
   of the affine subspace $\{ x = T_\nu \}\cap {\rm span}\, L$, such that $ {\cal B}'\cap {\cal L} = \varnothing$ and the point
$T_\nu \xi^{(\nu )}$ lies on the boundary of ${\cal B}'$. Let
$(T_\nu,\Xi_1,...,\Xi_n)$ be the center of ${\cal B}'$. Put
$$
\xi^{(\nu+1)} = \left( 1, \frac{\Xi_1}{T_\nu},...,
\frac{\Xi_n}{T_\nu}\right)\in {\rm span}\, L.$$ As ${\cal C}^{(\nu
+1)} \subset {\cal C}'$, we see that
 there are no nonzero points of $L$
in ${\cal C}^{(\nu+1)}$
 and (i) is valid with $\nu$ replaced by $\nu+1$.
 From the construction we see that (ii)  is also valid for  $\nu+1$. Lemma 4 is proved.

{\bf Remark 1.}\,\,\, Lemma 4 is obtained by well-known  arguments
(see \cite{SCH}, Chapter 3, \S 2). The constant $3^{-l-2}
(2w_{l-1}l !)^{-\frac{1}{l-1}} s^{\frac{1}{l-1}}$ in (\ref{prime})
may be slightly improved but this is   of no importance for the
proof  of our main result.

 \section{Lemmas concerning two sublattices}

{\bf Lemma 5.}\,\,\,{\it Let $\Gamma \subset \mathbb{Z}^{n+1}$ be
a  complete lattice, such that $  {\rm dim } \, \Gamma   = k+1\ge
3. $
 Let $R$ be
the $(k+1)$-dimensional fundamental volume of   $\Gamma $. Let
vector $\xi = (1,\xi_1,...,\xi_n) \in {\rm span }\,\Gamma $,
$\xi_j \in (0,1)$, be
  $(\Gamma,\gamma ,W)$-BAD. Consider a positive number $\kappa $, such that
\begin{equation}
\kappa \le \gamma W^{-\frac{k+1}{k}}.
 \label{kappasmall}
\end{equation}

Then there is a sublattice $\Lambda\subset \Gamma$, ${\rm dim }\,
\Lambda  = k$ satisfying the following two conditions:

1)
 the Euclidean distance from  $\xi \in \mathbb{R}^{n+1}$
 to ${\rm span }\, \Lambda \cap \{ x= 1\} $  does not exceed $\kappa$;

2) the $k$-dimensional fundamental volume $r$ of   $\Lambda$
admits the following upper bound:
 $$
r\le G (\gamma, \Gamma )\kappa^{-\frac{1}{k+1}},
$$
where \begin{equation}
 G(\gamma,\Gamma )=
   2^{2k^2 } w_k^{-k}w_{k-1} k!\gamma^{-\frac{k^3}{k+1}} R^k
 .
\label{je}
\end{equation}

}

Proof. We take
$$
M =(\gamma \kappa^{-1})^{\frac{k}{k+1}},\,\,\, \delta = M\kappa.
$$
By (\ref{kappasmall}) we have
$$
(M\gamma \delta^{-1})^{\frac{k}{k+1}} =(\gamma
\kappa^{-1})^{\frac{k}{k+1}} \ge W.
$$
We now can apply the statement 2) of Lemma 3 for $ l = k+1, L =
\Gamma$. Then (\ref{C1CB})  gives the inequality
\begin{equation}
\mu = \mu_k ({\cal C}_\xi (M, \delta))\le   \mu_{k+1} ({\cal
C}_\xi (M, \delta))\le 2^{2k}w_k^{-1} R\gamma^{-k}. \label{u}
\end{equation}
We see now that the cylinder $\mu\overline{\cal C}_\xi (M,
\delta)$ has $k$   linearly independent integer points
$\zeta_1,...,\zeta_k$. Define $\Lambda = {\rm span}
(\zeta_1,...,\zeta_k) \cap
 \mathbb{Z}^{n+1}$. The remark after Lemma 3 shows that the
 condition 1) is satisfied.
Let us obtain the needed upper bound for  the
 fundamental volume $r$ of $\Lambda $. As
 $$
 {\rm conv}(0,\zeta_1,..,\zeta_k) \subset {\rm span}\Lambda \bigcap
 \mu\overline{\cal C}_\xi (M,
\delta)
$$
(here ${\rm conv}(0,\zeta_1,..,\zeta_k)$ stands for the convex
hull of the points $0,\zeta_1,..,\zeta_k\in \mathbb{R}^{n+1}$)
 we see that
\begin{equation}
r\le k! \mu^k w_{k-1} \delta^{k-1}M, \label{v}
\end{equation}
 and the
required upper bound follows from (\ref{u},\ref{v}). Lemma is
proved.

  Let $\Gamma $ be a
sublattice  as in Lemma 5. For each lattice $\Gamma$  and for
every$\gamma$ small enough there exist $ (\Gamma, \gamma, W)$-BAD
vectors (see Example 1). Consider a $ (\Gamma, \gamma, W)$-BAD
vector $ \xi = (1,\xi_1,...,\xi _n) \in {\rm span } \,\Gamma$.
Then for any $ T\ge W$ the cylinder $ {\cal C}_\xi (T, \gamma
T^{-1/k}) $ contains no nonzero points of $\Gamma$. Let ${\cal B}$
be the facet $ \{ x = T\}$ of the cylinder $ {\cal C}_\xi (T,
\gamma T^{-1/k})\cap {\rm span } \,\Gamma$. This facet is a $k$
dimensional ball of radius $\gamma T^{-1/k}$
 centered at  $T\xi$.
Lemma 5 with
$$
\kappa = \frac{\gamma}{4n} \, T^{-\frac{k+1}{k}}
$$
implies that there is a $k$-dimensional complete sublattice
$\Lambda \subset \Gamma$ with $k$-dimensional fundamental volume
$r$ satisfying the condition
$$
 r\le G (\gamma, \Gamma )(4n\gamma^{-1})^{\frac{1}{k+1}} T^{\frac{1}{k}},
 $$
  and such that the intersection of  $ {\rm span } \,\Lambda$ with ${\cal B}$ is a  $(k-1)$-dimensional ball ${\cal B} '\subset {\cal B}$ with the
center $\Xi'$ and radius $\ge\gamma T^{-1/k}/2$.  Take another
$k-1$-dimensional ball ${\cal B} ''\subset {\cal B}'$  of radius
$2^{-3}\gamma T^{-1/k}$ centered at $\Xi '$. Then  the distance
from ${\cal B} ''$ to the boundary of  ${\cal B}$ is greater than
$2^{-3}\gamma T^{-1/k}$. Put $ \xi'=\Xi'/T$.
  Then the cylinder
\begin{equation}\label{lllw} {\cal C}_{\xi'} (T, 2^{-3}\gamma T^{-1/k} )\cap {\rm span }
\,\Lambda={\cal C}_{\xi'} (T, \gamma ' T^{-1/(k-1)} )\cap {\rm
span } \,\Lambda,\,\,\, \gamma' = 2^{-3}\gamma
T^{\frac{1}{k(k-1)}}
\end{equation}
 contains no nonzero points of  $\Lambda$ and
we can apply Lemma 2 with $ l = k $ and $L=\Lambda$. Thus we
obtain a $(\Lambda,\hat{\gamma }, T)$-BAD vector $\hat{\xi }$ with
\begin{equation}
\hat{\gamma } = \hat{\gamma } ( \gamma , T, \Lambda)
 = \gamma^* (2^{-3}\gamma
T^{\frac{1}{k(k-1)}},\Lambda ). \label{hat}
\end{equation}
As the cylinder (\ref{lllw}) has no nonzero points of $\Lambda$ we
see from the Minkowski theorem on convex bodies that
$$
2T\cdot w_{k-1} \left(2^{-3} \gamma T^{-\frac{1}{k}}\right)^{k-1}
< 2^k r
$$
or
\begin{equation}
r > 2^{-2(k-1)} T^{\frac{1}{k} }w_{k-1} \gamma^{k-1}.  \label{vo}
\end{equation}

Note that from (\ref{hat},\ref{prime},\ref{vo}) we see that
\begin{equation}
\hat{\gamma }  \ge \underline{C}(k)\gamma
T^{\frac{1}{k(k-1)}},\label{hatmin}
\end{equation}
where
\begin{equation}
\underline{C}(k)= \min\left(3^{-5} , 2^{-1-\frac{1}{k-1}} 3^{-k-1}
(k!)^{-\frac{1}{k-1}}\right).
 \label{ce}
\end{equation}
Now we put
\begin{equation}
 Z_1 (\gamma, k ) = (\underline{C}(k)\gamma )^{-\frac{n}{n+1}},
  \label{zet1}
\end{equation}
\begin{equation}
   Z_2 (i,\gamma, \Gamma ) =   \left( i\cdot 2^{2(k-1)}w_{k-1}^{-1}
G(\gamma,\Gamma ) (4n\gamma^{-1})^{\frac{1}{k+1}}
(\underline{C}(k)\gamma )^{-\frac{(k-1)^2}{k}}
\right)^{\frac{nk}{n+1-k}}.
 \label{zet2}
\end{equation}

{\bf Lemma 6.}\,\,\,{\it For the vector $\hat{\xi}$ defined above
and for  $N $ under the condition
\begin{equation}\label{starr}
N \ge Z(i,\gamma,\Gamma, T) = \max (  Z_1 (\gamma, k ) T^{\frac{n(k+1)}{(n+1)k}},
  Z_2 (i,\gamma, \Gamma ) T^{\frac{n }{  k(n+1-k)}})
\end{equation}
we have
$$
\mu_k (\hat{\xi}, N) \le i^{-1}.
$$}

Proof. As $ N \ge  Z_1 (\gamma, k ) T^{\frac{n(k+1)}{(n+1)k}}$, we
see from (\ref{zet1},\ref{ce},\ref{hatmin}) that the condition
(\ref{EN}) of the case 2) of Corollary 1 to Lemma 3 is  satisfied.
Then due to (\ref{C1C})  we have
$$
\mu_k (\hat{\xi}, N) \le 2^{2(k-1)} w_{k-1}^{-1} r \hat{\gamma
}^{-\frac{(k-1)^2}{k}} N^{\frac{k-n-1}{nk}}.
$$
It remains to make use of the inequality  $ N \ge
  Z_2 (i,\gamma, \Gamma ) T^{\frac{n }{  k(n+1-k)}}
$ and of the formulas (\ref{zet2},\ref{vo},\ref{hatmin},\ref{ce}).
Lemma 6 follows.

We shall need the following notation related to a pair of
sublattices.

Let $\Lambda \subset \Gamma \in \mathbb{Z}^{n+1}$ be complete
sublattices such that

$$
  {\rm dim}\, \Gamma
 >{\rm dim}\,\Lambda  .
$$
Then $\Gamma$ can be partitioned into classes $ \pmod{\Lambda}$:
$$
\Gamma = \bigcup_{\alpha \in \mathbb{Z}^v} \Gamma_\alpha,\,\,\,
\Gamma_0 =\Lambda,\,\,\, v = {\rm dim}\,(\,{\rm span }\,\Gamma )
 -{\rm dim}\,(\,{\rm span }\,\Lambda  ),
$$
so that the affine subspaces ${\rm aff}\, \Gamma_\alpha$ are
parallel   ${\rm span }\,\Lambda ={\rm aff}\, \Gamma_0$. Here by
${\rm aff}\, \Omega$ we mean the smallest affine subspace of
$\mathbb{R}^{n+1}$ containing $\Omega$.

 Denote by $R = R (\Lambda,\Gamma ) >0$ the minimal distance
between points $ z^{(1)},z^{(2)}$, where $ z^{(1)}\in
\Gamma\setminus \Lambda $ and $z^{(2)} \in {\rm span }\,\Lambda$.
For our purpose we need not the  $R (\Lambda, \Gamma )$ itself but
a little bit different distance $\rho = \rho (\Lambda,\Gamma )$
which we define now. Recall that as it was pointed out in the very
beginning of Section 3 all $k$-dimensional lattices $\Lambda$
under the consideration admit the property
$$
{\rm dim}\,( \Lambda \cap {\cal P} ) = k-1
$$
(here ${\cal P} =\{ z=(x,y_1,...,y_n)\in \mathbb{R}^{n+1}:\,\, x =
1\}$). For such  a lattice $\Lambda$  and for a lattice $\Gamma
\supset \Lambda$  of dimension $ {\rm dim}\, \Gamma > {\rm dim}\,
\Lambda$ we consider the following objects.
 Put  ${\cal L} = {\rm span
}\,\Lambda \cap {\cal P}$. Let ${\cal G}$ be the parallel
projection of   $\Gamma \setminus \Lambda$ along $ {\rm span
}\,\Lambda $ onto ${\cal P}$. Then we define $\rho = \rho
(\Lambda,\Gamma  )$ to be the minimal distance between points $
z^{(1)},z^{(2)}$, where $ z^{(1)}\in {\cal G}$ and $z^{(2)} \in
{\cal L} = {\rm span }\,\Lambda \cap {\cal P}$. We see that $\rho
= \rho (\Lambda,\Gamma  )  >0$. In fact, $ \rho (\Lambda,\Gamma
)\ge \rho (\Lambda , \mathbb{Z}^{n+1} )
>0$.

Now we give two more lemmas.

{\bf Lemma 7.}\,\,\,{\it Let $\Lambda \subset \Gamma\subset
\mathbb{Z}^{n+1}$ be complete sublattices, such that
$$
   {\rm dim} \, \Gamma  = k+1 = {\rm dim}\,   \Lambda
+1 ,$$ and
 $\rho = \rho (\Lambda,\Gamma  ).
 $
Let $\xi = (1,\xi_1,...,\xi_n) \in {\rm span }\,\Lambda$ be a
$(\Lambda ,\gamma ,W )$-BAD vector with some positive $\gamma$ and
$W\ge 1$.    Put
\begin{equation}
\gamma ' =\gamma ' (\gamma , \Lambda , \Gamma )=
\gamma^{\frac{k-1}{k}} 2^{-\frac{1}{k}} \rho^{+\frac{1}{k}}
\label{gammaprime}
\end{equation}
and
\begin{equation}
  A_1 = A_1(\gamma , \Lambda , \Gamma )
= \max\left(   (\rho (2\gamma ')^{-1})^{\frac{k}{k-1}} , (2\gamma
' \rho^{-1})^k\right). \label{AAA}
\end{equation}
 Suppose that
\begin{equation}
 T \ge   A_1W^{\frac{k}{k-1}}
\ge\max\left( (\rho (2\gamma ')^{-1}W)^{\frac{k}{k-1}} , (2\gamma
' \rho^{-1})^k\right).
  \label{te}
\end{equation}
Let $\xi' = (1,\xi_1',...,\xi_n') \in {\rm span }\,\Gamma $
satisfy the following two conditions:

1) the orthogonal projection of vector $\xi '$ on the subspace
${\rm span }\,\Lambda$ is of the form $\lambda \xi $ with some
positive $\lambda$;

2) for the Euclidean norm we have $| \xi ' -  \xi | =
(2T)^{-1}\rho$.

Then
$$
{\cal C}_{\xi '} (T,\gamma' T^{-1/k}) \bigcap \Gamma = \{ 0\}.
$$}

 Proof.
The $(k+1$-dimensional linear subspace ${\rm span}\, \Gamma$
contains parallel $k$-dimensional affine subspaces $ {\rm aff}\,
\Gamma_i, i\in \mathbb{Z}$. Each such subspace $ {\rm aff}\,
\Gamma_i$  divides the subspace ${\rm span}\, \Gamma$ into two
"half-subspaces" with the common boundary $ {\rm aff} \,\Gamma_i$.
The situation with the $k$-dimensional affine subspace ${\rm
span}\, \Gamma \cap \{ x = T\}$ and $(k-1)$-dimensional affine
subspaces $ {\rm aff}\, \Gamma_i \cap \{ x = T\}, i\in \mathbb{Z}$
is quite similar. We should note that the Euclidean distance
between neighboring subspaces
 $ {\rm aff}\, \Gamma_i \cap \{ x = T\}$ and  $ {\rm aff}\, \Gamma_{i+1} \cap \{ x =
 T\}$ is exactly $\rho$.

 Without loss of generality we may suppose that the point
  $T\xi'\in {\rm span }\,\Gamma \cap \{ x = T\}$ lies
in the same "half-subspace" (with respect to ${\rm aff }\,\Gamma_0
\cap \{ x = T\} = {\rm span }\,\Lambda \cap \{ x = T\}$) as the
set ${\rm span }\,\Gamma_1 \cap \{ x = T\}$. From the conditions
1), 2) we see that the distance from  the point $T\xi'$ to the
subspace ${\rm aff }\,\Gamma_0 \cap \{ x = T\}$ is equal to the
distance from $T\xi'$ to the subspace ${\rm aff }\,\Gamma_1 \cap
\{ x = T\}$ and is equal to $\rho/2$.

  Define $ H = 2\gamma'\rho^{-1} T^{\frac{k-1}{k}} $. Then
by definition of $H$ and (\ref{te}), we have $H\le T$. Note that
the distance from each point of the form $t\xi, H\le t\le T$ to
the corresponding affine subspace ${\rm span }\,\Lambda \cap \{
x=t\}$ is greater than $\gamma' T^{-1/k}$. So
$$
{\cal C}_{\xi '} (T,\gamma' T^{-1/k})\bigcap \{ z=(x,y_1,...,y_n)
: |x| \ge H\} \bigcap {\rm span}\, \Lambda = \varnothing.
$$
It means that the cylinder ${\cal C}_{\xi '} (T,\gamma' T^{-1/k})$
intersected with the domain $\{ H\le x\le T\}$ has no points of
the lattice $\Lambda$. But the distance from each point of the
form $t\xi, H\le t\le T$ to the corresponding affine subspace
${\rm aff }\,\Gamma_1 \cap \{ x=t\}$ is greater than the distance
from $t\xi $ to   ${\rm span }\,\Lambda \cap \{ x=t\}$. So the
distance from each point of the form $t\xi, H\le t\le T$ to the
any  affine subspace ${\rm aff }\,\Gamma_i \cap \{ x=t\}, i \neq
0$
 is greater than $\gamma' T^{-1/k}$ also. So
$$
{\cal C}_{\xi '} (T,\gamma' T^{-1/k})\bigcap \{ z=(x,y_1,...,y_n)
: |x| \ge H\} \bigcap \Gamma = \varnothing.
$$
From the other hand if $ 0\le t\le H$ then the distance from
 $t\xi$ to any ${\rm aff }\,\Gamma_i \cap \{ x=t\}, i \neq
0$ is again greater than  $\gamma' T^{-1/k}$. Hence
$$
{\cal C}_{\xi '} (T,\gamma' T^{-1/k})\bigcap \Gamma = {\cal
C}_{\xi '} (H,\gamma' T^{-1/k})\bigcap \Lambda.
$$
But  (\ref{gammaprime}) implies that $ \gamma ' T^{-1/k} =\gamma
H^{-1/(k-1)}$. As $\xi$ is a $(\Lambda ,\gamma, W )$-BAD vector
we see that
$$
{\cal C}_{\xi '} (H,\gamma' T^{-1/k})\bigcap \Lambda\subseteq
{\cal C}_{\xi } (H,\gamma H^{-1/(k-1)})\bigcap \Lambda = \{ 0\}.
$$
(Note that from (\ref{te}) it follows that $ H\ge W$.)
  Lemma 7 is proved.

 {\bf Lemma 8.}\,\,\,{\it In the notation of Lemma 5, let $r$ be
the $k$-dimensional fundamental volume of the lattice $\Lambda $,
vector $\xi '$ be defined in Lemma 7 and let $\xi''=
(1,\xi_1'',...,\xi_n'') \in {\rm span }\,\Gamma $ be a vector
satisfying
\begin{equation}
|\xi''-\xi ' |< \rho (4T)^{-1}. \label{ex}
\end{equation}
Set
\begin{equation}
A_2 =A_2 (\gamma, \Lambda ,\Gamma ) = 3\rho \gamma^{-1}/4,
\label{aaa}
\end{equation}
\begin{equation}
B_1 = B_1  (\gamma )=
 (\sqrt{2}\gamma )^{-\frac{n}{n+1}},\,\,
B_2 = B_2(\Lambda ,\Gamma )=
\left(\frac{2\sqrt{2}}{3\rho}\right)^{\frac{n}{n+1}}, \label{BBBB}
\end{equation}
\begin{equation}
C_1 = C_1(\gamma, \Lambda ) = 2^{2k-2}w_{k-1}^{-1}
r\gamma^{1-k},\,\,
 C_2= C_2  (\gamma, \Lambda  )
=2^{\frac{4k^2-3k-1}{2k}}w_{k-1}^{-1}r\gamma^{-\frac{(k-1)^2}{k}},
\label{ceodin}
\end{equation}
$$
  C_3= C_3  (\gamma, \Lambda ,\Gamma )
=2^{\frac{k^2-3k-1}{2k}} 3^ {\frac{1}{k}}
w_{k-1}^{-1}r\gamma^{-\frac{(k-1)^2}{k}}\rho^{\frac{1}{k}}.
$$

 Suppose that
\begin{equation}
 T \ge   A_2W^{\frac{k}{k-1}}.
   \label{te0}
\end{equation}
Then the following statements are valid:

1) for $N$ in the interval
\begin{equation}
N \le B_1 W^{\frac{kn}{(k-1)(n+1)}}
     \label{en00}
\end{equation}
we have
\begin{equation}
\mu_{k} (N, \xi '') \le  C_1 W N^{-1};
 \label{LL00} \end{equation}

2) for $N$ in the interval
\begin{equation}
B_1 W^{\frac{kn}{(k-1)(n+1)}}\le
   N\le B_2T^{\frac{n}{n+1}}  \label{en0}
\end{equation}
we have
\begin{equation}
\mu_{k} (N, \xi '') \le  C_2 N^{\frac{k-n-1}{nk}}  ;
 \label{LL0} \end{equation}

3) for $N$ in  the interval
\begin{equation}
   N\ge B_2T^{\frac{n}{n+1}} \label{en1}
\end{equation}
we have
\begin{equation}
\mu_{k} (N, \xi '') \le  C_3 T^{-\frac{1}{k}} N^{\frac{1}{n}}.
 \label{LL} \end{equation}

 }

{\bf Corollary.}\,\,\,{\it  Under the  conditions of Lemma 8, for
$N$ in the interval
\begin{equation}
H(i,\gamma , \Lambda, W) =  \max\left( (C_1 (\gamma, \Lambda ) i
W) , (C_2 (\gamma , \Lambda ) i)^{\frac{nk}{n+1-k}}\right) \le
N\le (i  C_3 (\gamma, \Lambda ,\Gamma ))^{-n} T^{\frac{n}{k}}
\label{inte}
\end{equation}
we  have the following inequality:
\begin{equation}
\mu_{k} (N, \xi '') \le   i^{-1}.
 \label{mumain} \end{equation}}

 Proof of Lemma 8. First of all, let us consider the case 3).

 Set
$$
M = B_2^{\frac{n+1}{n}}TN^{-\frac{1}{n}}\le N,\,\,\,\delta
=N^{-\frac{1}{n}}/\sqrt{2}.
$$
It follows from the definition of $\xi '$ and (\ref{ex})  that
\begin{equation}
 {\cal C}_{\xi ''} (N, N^{-1/n}) \supset {\cal C}_{\xi  } (M,
\delta)\cap {\rm span }\,\Lambda . \label{false}
\end{equation}
  Hence
$$
\mu_{k} (N, \xi '')\le \mu_{k}({\cal C}_{\xi  } (M, \delta)),
$$
so it is suffices to obtain the corresponding upper bound for the
latter successive minimum. We observe that by (\ref{te0}) we have
$$
(M\delta^{-1}\gamma)^{\frac{k-1}{k}} = \left(\frac{4}{3}
T\rho^{-1}\gamma \right)^{\frac{k-1}{k}}\ge W.
$$
Applying the statement 2) of Lemma 3 we obtain (\ref{C1CB}) with $
l=k, s=r$. Now (\ref{LL}) follows from (\ref{C1CB}).

Consider the case 2). Since $ M \ge N$, the relation (\ref{false})
may be false, so  we have
\begin{equation}
 {\cal C}_{\xi ''} (N, N^{-1/n}) \supset {\cal C}_{\xi  } (N,
\delta)\cap {\rm span }\,\Lambda . \label{falseno}
\end{equation}
  Hence
$$
\mu_{k} (N, \xi '')\le \mu_{k}({\cal C}_{\xi  } (N, \delta)),
$$
It follows from (\ref{en0})  that
$$
(N\delta^{-1}\gamma)^{\frac{k-1}{k}} \ge W.
$$
Let us apply the statement 2) of Lemma 3 for the cylinder ${\cal
C}_{\xi  } (N, \delta)$ from (\ref{falseno}). Then the conclusion
(\ref{C1CB}) of Lemma 3 with our parameters leads to (\ref{LL0}).

Finally, we consider the case 1). Again, we have $M\ge N$. So we
must use the relation (\ref{falseno}). But  (\ref{en00}) implies
that
$$
(N\delta^{-1}\gamma)^{\frac{k-1}{k}} \le W.
$$
Applying the statement 1) of Lemma 3 we get from  (\ref{C1C0B})
the desired inequality (\ref{LL00}).

Lemma 8 is proved.

\section{Proof of Theorem 1: general case}

  Now we are able to  give a proof of Theorem 1 in the case $ k\ge 2$. We begin  with the same consideration of the countable set
 of all the $n$-dimensional complete sublattices of the integer lattice
$\mathbb{Z}^{n+1}$. We fix an enumeration of this set and let
$$ L_1, L_2,..., L_i , ...
$$
be all these lattices. Set $\pi_i = {\rm span} \,L_i$. Suppose
that $$\pi_1 = \{z=(x,y_1,...,y_n)\in\mathbb{R}^{n+1}:\,\ x=0\}.$$

Let $2\le k \le n-1$. We construct a  sequence of real numbers
$$
\eta_1> \eta_2>...>\eta_i >...
$$
decreasing to zero, a sequence of positive real numbers
$$
\gamma_1,\gamma_2,...,\gamma_i,... ,
$$
two sequences of real numbers
$$
W_1,W_2,...,W_i,...,
$$
$$
H_1,H_2,...,H_i, ... ,$$
$$W_i    \ge 1,\,\,\,
W_i, H_i \to +\infty,\,\,\, i \to +\infty,
$$
 two sequences
 of complete sublattices
$$
\Lambda_1,\Lambda_2,...,\Lambda_{i-1},\Lambda_i , ... ,
$$
$$
\Gamma_2,\Gamma_3,...,\Gamma_{i},\Gamma_{i+1},... ,
$$
and a sequence of vectors
$$
\xi_i = (1,\xi_{i,1},..., \xi_{i,n})\in \mathbb{R}^{n+1}
$$
satisfying the following conditions (i) -- (vii). Further, suppose
$r_i$  be the $k$-dimensional fundamental volume of $\Lambda_i$
and let $R_i$ be the $(k+1)$-dimensional fundamental volume of
$\Gamma_i$.

(i)  For every $i\in \mathbb{N}$ we have
$$
\Lambda_i \subset \mathbb{Z}^{n+1}, \,\,\, {\dim }\, \Lambda_i =
k;
$$
$$
\Gamma_{i+1} \subset \mathbb{Z}^{n+1},\,\,\, {\dim }\,
\Gamma_{i+1} = k+1;
$$
$$
\Lambda_{i}, \Lambda_{i+1} \subset \Gamma_{i+1}.
$$

(ii) For every $i\in \mathbb{N}$ the vector $\xi_i$ is $(\Lambda_i
,\gamma_i, W_i)$-BAD.

(iii) The $n$-dimensional closed ball $\overline{\cal B}_i\subset
\{ z = (x,y_1,..., y_n) \in \mathbb{R}^{n+1}:\,\, x= 1\}$ of
radius $\eta_i$ is centered at $\xi_i$ and   has no common points
with   $\pi_{i}$.

(iv) The balls defined in (iii) form a nested sequence
$$
\overline{\cal B}_1 \supset \overline{\cal B}_2 \supset ...
\supset \overline{\cal B}_i .
$$

(v) For every $i\ge 2$ the following inequality holds:
\begin{equation}
 H_i  \ge
\max \left( H(2(i+1),\gamma_i,\Lambda_i, W_i), \frac{4(i+1)}{\rho
(\Lambda_i, \mathbb{Z}^{n+1} )}\right)
 \label{AH10}
\end{equation}
(here the value of $\rho (\cdot, \cdot)$ for two lattices is
defined in Section 4 before Lemma 7 and $H(\cdot,\cdot,\cdot, \cdot)$
is defined in (\ref{inte})).

 (vi) For every $i\ge 2$, every  $\xi \in
\overline{\cal B}_i$ and for every real $N$ in the interval
$H_{i-1} \le N< H_i $
    one
has
$$
\mu_k (\xi, N) \le i ^{-1} .
$$

(vii) For every $i\ge 2$, every $\xi \in \overline{\cal B}_i$ and
every real $N$ in the interval $H_{i-1}^n \le N< H_i^n $
    one
has
$$
\mu_{k+2} (\xi, N) \ge i .
$$

  Suppose that
all these objects are already  constructed. Then we have  Theorem
1  proved in the case $ k \ge 2$. Indeed, if we consider the
unique vector $\xi =(1,\xi_1,...,\xi_n)$ from the intersection $
\cap_{i\in \mathbb{N}} \overline{\cal B}_i$, then the components
   $1,\xi_1,...,\xi_n $ are linearly
independent over $\mathbb{Z}$ due to (iii), and
$$\lim_{N \to +\infty }
\mu_k (\xi, N) =0,\,\,\,\,\lim_{N \to +\infty } \mu_{k+2} (\xi, N)
=+\infty
$$
due to (vi) and  (vii).

We now  describe an inductive process, which constructs all  the
objects mentioned.

First of all, put $ W_1=H_1=1$,
$$
\Lambda_1 = \mathbb{Z}^{n+1}\bigcap
\{z=(x,y_1,...,y_n)\in\mathbb{R}^{n+1}:\,\,\, y_k = ...=y_n = 0\},
$$
 Take $\xi$ to be a $(\Lambda,\gamma_1,1)$-BAD vector with some positive $\gamma_1$
(we can take such a vector from Example 1). We do not define
  $\Gamma_1$.

Obviously, $\rho (\Lambda_1, \mathbb{Z}^{n+1}) = 1 $.
 The conditions (i) -- (vii) for $ i = 1$ are satisfied (note that the conditions (v) -- (vii) are empty).

Assume that all the objects $ \eta_i , \gamma_i,  W_i, H_i
,\Lambda_i , \Gamma_{i}, \xi_i $
 for every natural $i$ up to $t$
    are constructed
to  satisfy  the conditions (i) -- (vii). Let us describe the
construction for $ i = t+1$.

Consider the integer vector $q\in \mathbb{Z}^n$ orthogonal to the
subspace $\pi_{t+1}$. For any small positive $\varepsilon$ there
exists
 an integer vector $q'$ such that the angle between $q$ and $q'$ is less than $\varepsilon$ and
$q' \not\in \Lambda_t$. Put
$$
\Gamma_{t+1} = \mathbb{Z}^{n+1} \cap {\rm span}\, (\Lambda_t\cup
q').
$$
Then $ \Gamma_{t+1}$ is a complete sublattice of dimension  ${\rm
dim}\, \Gamma_{t+1} = k+1$, $\Gamma_{t+1} \supset \Lambda_t$,  and
  $$ {\rm span} \,\Gamma_{t+1} \not\subset \pi_{t+1}.
$$
Moreover for the vector $ q' \in \Gamma_{t+1}$ and any nonzero
vector $p\in \pi_{t+1}$ the angle between $ p$ and $q'$ is greater
than $\frac{\pi}{2}-\varepsilon$.

Let $R_{t+1}$ be the $(k+1)$-dimensional fundamental volume of
$\Gamma_{t+1}$. Set
$$
\rho_t^{(1)} = \rho (\Lambda_t, \mathbb{Z}^{n+1}),\,\,\,
\rho_t^{(2)} = \rho (\Lambda_t, \Gamma_{t+1}),\,\,\, \rho_t^{(3)}
= \rho (\Gamma_{t+1}, \mathbb{Z}^{n+1}).
$$
Set
$$
E_j(t) = A_j(\gamma_t ,
\Lambda_t,\Gamma_{t+1})W_t^{\frac{k}{k-1}},\,\, j =1,2,
$$
where the right hand sides are defined by (\ref{AAA},\ref{aaa}),
and set
$$
E_3(t) = \frac{3\rho^{(2)}_t}{4\eta_t},\,\,\, E_4(t) =
\frac{2^{n+3}(t+1)^{n+1} \rho^{(2)}_t}{\rho^{(1)}_t
(\rho^{(3)}_t)^n}.
$$
We   also need one more quantity $E_5(t)$ defined as follows.
First, we put
$$
Z_1(t) = Z_1 (\gamma^* (\gamma_t^{\frac{k-1}{k}} 2^{-\frac{1}{k}}
(\rho_t^{(2)})^{\frac{1}{k}}, \Gamma_{t+1}),k),
$$
$$
Z_2(t) = Z_2 (2(t+1), \gamma^* (\gamma_t^{\frac{k-1}{k}}
2^{-\frac{1}{k}} (\rho_t^{(2)})^{\frac{1}{k}},
\Gamma_{t+1}),\Gamma_{t+1}),
$$
where $Z_1(\cdot,\cdot ),Z_2(\cdot,\cdot,\cdot )$ are defined by
(\ref{zet1},\ref{zet2}) and $\gamma^*(\cdot,\cdot  )$ is defined
by (\ref{prime}). Then we put
$$
E_5(t)= \max \left( (Z_1(t))^{\frac{(n+1)k}{n(n-k)}} (2(t+1) C_3
(\gamma_t,\Lambda_t,\Gamma_{t+1}))^{\frac{(n+1)k}{n-k}} ,
(Z_2(t))^{\frac{(n+1-k)k}{n(n-k)}} (2(t+1) C_3
(\gamma_t,\Lambda_t,\Gamma_{t+1}))^{\frac{(n+1-k)k}{n-k}}\right),
$$
where $C_3(\cdot,\cdot,\cdot )$ is defined by (\ref{ceodin}). Note
that as $k\le n-1$ all the exponents are positive (particulary, $
n-k \ge 1$ and all the denominators in the exponents are nonzero).

 Put
$$
T_t = \max_{1\le j \le 5} E_j (t).
$$
Since $ T_t \ge E_1(t), E_2 (t)$, we can apply Lemmas 7,8 to the
lattices $\Lambda = \Lambda_t, \Gamma = \Gamma_{t+1}$. Denote by
$\xi_t '$ the real $(n+1)$-dimensional vector satisfying the
conditions 1), 2) of Lemma 7. Consider the ball
$$
{\cal B}_t' = \{ \xi = (1,\xi_1,...,\xi_n) :\,\, |\xi - \xi_t' | <
\rho^{(2)}_t (4T_t)^{-1}\}.
$$
Note that for $\xi \in {\cal B}_t' $ we have $|\xi - \xi_t| <
\frac{3\rho^{(2)}_t}{4T_t}.$
 Since $ T_t \ge E_3(t)$, we have
$$
{\cal B}_t' \subset  {\cal B}_t .
$$
Note that
$$
H_t \ge H(2(t+1),\gamma_t,\Lambda_t,W_t)
$$
by the inductive conjecture (v).  So by  Corollary of Lemma 8 we see that for every $N$ in  the interval
$ H_t \le N \le H_t',\,\,\, $ where
\begin{equation} H_t' = \left(2 (t+1) C_3 (\gamma_t, \Lambda_t, \Gamma_{t+1})\right)^{-n }T_t^{\frac{n}{k}},
\label{ahtprim}
\end{equation}
and every $\xi \in {\cal B}_t'\cap {\rm span}\,\Gamma_{t+1}$ we
have
$$
\mu_k (\xi, N) \le (2(t+1))^{-1}.
$$

Let us prove that for any $N\ge H_t^n$ and for any $\xi \in {\cal
B}_t'\cap {\rm span} \,\Gamma_{t+1} $ we have
\begin{equation}
\mu_{k+2} (\xi , N) \ge 2(t+1). \label{dd}
\end{equation}
To do this let us put
$$
U = T_t \cdot \frac{\rho^{(1)}_t}{3(t+1)\rho^{(2)}_t} .
$$
Recall that
$$
|\xi_t - \xi_t'| =\frac{\rho_t^{(2)}}{2T_t}
$$
and $\xi_t \in {\rm span }\,\Lambda_t$. So for every $N$ in the interval  $ H_t^n \le N \le U$ we have
$$
|N\xi - N\xi_t| \le U \cdot \frac{3\rho^{(2)}_t}{4T_t}\le
\frac{\rho^{(1)}_t}{4(t+1)} ,
$$
and so the distance between $N\xi $ and  $ {\rm span } \Lambda $
does not exceed $ \frac{\rho^{(1)}_t}{4(t+1)} $. But it follows
from the condition (v) of the $t$-th step that for the considered
values of $N$ we have
$$
N^{-\frac{1}{n}} \le H_t^{-1} \le \frac{\rho^{(1)}_t}{4(t+1)}.
$$
We see that the maximal distance  between a point of the cylinder
 $ C_\xi (N, N^{-1/n})$ and the subspace ${\rm span}\, \Lambda_t$
 is $\le \frac{\rho^{(1)}_t}{2(t+1)}$.
 Recall that ${\rm dim}\, \Lambda_t = k$.
Hence the cylinder $ 2(t+1)\cdot C_\xi (N, N^{-1/n})$ cannot
contain $ k+1$ linearly independent integer points inside  for $
H_t^n \le N \le U$, so in this case we have the inequality
$$
\mu_{k+1} (\xi , N) \ge 2(t+1)
$$
(and thus, the inequality (\ref{dd})).

Suppose that $ N \ge U$. Then    we deduce from the inequality $
T_t \ge E_4(t)$ that
$$
N^{-\frac{1}{n}} \le U^{-\frac{1}{n}} \le
\frac{\rho^{(3)}_t}{2(t+1)}.
$$
So the distance between a point of
 $ C_\xi (N, N^{-1/n})$ and the linear subspace ${\rm span}\,
 \Gamma_{t+1}$ is
 $\le \frac{\rho^{(3)}_t}{2(t+1)}$.
 Hence the cylinder
 $ 2(t+1)\cdot C_\xi (N, N^{-1/n})$ cannot
contain $ k+2$ linearly independent integer points inside. This
implies  (\ref{dd})   in the case $ N \ge U$.

 We have proved the following statement: for any $ \xi \in {\cal B}_t'\cap {\rm span}\, \Gamma_{t+1}$ we have
\begin{equation}
\mu_{k+2} (\xi , N) \ge 2(t+1),\,\,\, N \ge H_t^n,
\label{kaplusdva}
\end{equation}
\begin{equation}
\mu_{k} (\xi , N) \le (2(t+1))^{-1},\,\,\,  H_t \le N \le H_t ',
\label{ka}
\end{equation}
where $ H_t'$ is defined by  (\ref{ahtprim}). Moreover, if
$$ \gamma '_t = \gamma '(\gamma_t , \Lambda_t , \Gamma_{t+1} ),
$$
where $\gamma  ' (\cdot , \cdot , \cdot )$ is defined by
(\ref{gammaprime}) then it follows from Lemma 7
that  the cylinder
$$
{\cal C}_{\xi_t '} (T_t, \gamma_t' T_t^{-1/k} )
$$
contains  no nonzero points of $\Gamma_{t+1}$. Recall that we
constructed   $\Gamma_{t+1}$ to satisfy the condition $ {\rm span}
\,\Gamma_{t+1} \not\subset \pi_{t+1}. $ Moreover for the vector $
q' \in \Gamma_{t+1}$ and any nonzero vector $p\in \pi_{t+1}$ the
angle between $ p$ and $q'$ is greater than
$\frac{\pi}{2}-\varepsilon$.
 We can find a
$k$-dimensional  ball ${\cal }B'' $ of   radius $\gamma_t'
T_t^{-1/k}/2$ inside the facet $\{ x = T_t\}\cap {\rm span}\,
\Gamma_{t+1}$ of the cylinder $ {\cal C}_{\xi_t '} (T_t, \gamma_t'
T_t^{-1/k} ) \cap {\rm span}\, \Gamma_{t+1} $ (in fact, this facet
is a $k$-dimensional  ball of radius $\gamma_t' T_t^{-1/k}$) such
that ${\cal }B''\cap ({\rm span}\, \Gamma_{t+1}\cap \pi_{t+1}\cap
\{ x = T_t\})=\varnothing$. Let $\Xi''$ be the center of ${\cal
}B'' $. Put $\xi_t'' = \Xi'' / T_t$. Then $ \xi_t''\in {\rm
span}\, \Gamma_{t+1}$. From the construction we see that
$n$-dimensional ball with the center at  point $\Xi''$ and radius
$\gamma_t' T_t^{-1/k}/4$ has no common points with the subspace
$\pi_{t+1}$.

We get a cylinder
$$
{\cal C}_{\xi_t ''} (T_t, \gamma_t' T_t^{-1/k}/4 )
$$
with
 no nonzero points of $\Gamma_{t+1}$
inside it and with    $\xi_t'' \in {\rm span}\, \Gamma_{t+1}$.
 By Lemma 4    we construct a $(\Gamma_{t+1},\gamma^*_t,
T_t)$-BAD vector $\xi^*_t \in {\rm span}\, \Gamma_{t+1}$ with
$$
\gamma_t^* = \gamma^* (\gamma_t'/4, \Gamma_{t+1} )
$$
( $\gamma ^* ( \cdot , \cdot )$ defined by (\ref{prime})), such
that the facet $\{ x = T_t\} $ of
$$
{\cal C}_{\xi_t^*} (T_t, \gamma_t^* T_t^{-1/k} )
$$
lies inside the facet $\{ x = T_t\} $ of   $ {\cal C}_{\xi_t '}
(T_t, \gamma_t' T_t^{1/k} ) $ and
  does not intersect  $\pi_{t+1}$. Hence the ball
 $$
 {\cal B}^*_t = \{ \xi = (1,\xi_1,...,\xi_n) \in \mathbb{R}^{n+1}:\,\,\,
 |\xi - \xi^*_t|< \gamma_t^* T_t^{-(k+1)/k}\}
 $$
 enjoys the following properties:
\begin{equation}
{\cal B}^*_t\subset {\cal B}_t',\,\,\, {\cal B}^*_t\bigcap
\pi_{t+1} = \varnothing . \label{ball}
\end{equation}

Recall that    $\xi^*_t \in {\rm span}\, \Gamma_{t+1}$ is  a $
(\Gamma_{t+1},\gamma^*_t, T_t)$-BAD vector.
 Applying Lemma 5
to the lattice $ \Gamma =\Gamma_{t+1}$, $
(\Gamma_{t+1},\gamma^*_t, T_t)$-BAD vector $\xi^*_t$ and $\kappa =
\gamma_t^* T_t^{-(k+1)/k}/(4n)$
 we get a complete $k$-dimensional lattice $\Lambda_{t+1}$ with fundamental volume
\begin{equation}
r_{t+1} \le  G (\gamma_t^*, \Gamma_{t+1} )
(\gamma_t^*/4n)^{-\frac{1}{k+1}} (T_t)^{\frac{1}{k}} \label{voll}
\end{equation}
(here $G(\cdot ,\cdot )$ is defined by (\ref{je})),
  such that
$ \Lambda_{t+1} \subset \Gamma_{t+1} , $
 and
 the Euclidean distance between   $\xi^*_t  $
and   ${\rm span }\, \Lambda_{t+1} \cap \{ x= 1\} $  does not
exceed $\gamma_t^* (T_t)^{-(k+1)/k}/4n$.

Next, we apply the construction described in Section 4 after Lemma
5 and
 obtain a
$(\Lambda_{t+1},\gamma_{t+1}, T_{t+1})$-BAD vector
$$
\xi_{t+1} \in {\rm span}\,\Lambda_{t+1}.
$$
We set
 $$W_{t+1} = T_t.
 $$

In the notation of  Section 4 we have
$$
\xi_{t+1} =\hat {\xi^*_t},
$$
$$
\gamma_{t+1} =  \hat{\gamma } ( \gamma_t^*, T_t, \Lambda_{t+1})=
\gamma^* (2^{-3} \gamma_t^*T_t^{\frac{1}{k(k-1)}}, \Lambda_{t+1})
$$ (here $
 \hat{\gamma } ( \cdot, \cdot,\cdot )$ is defined by (\ref{hat})
and $\gamma ^* ( \cdot,\cdot )$ is defined by (\ref{prime})).

Now we set
\begin{equation}
H_{t+1} = \max \left( Z(2(t+2),\gamma_t^*,\Gamma_{t+1}, T_t),
H(2(t+1),\gamma_{t+1},\Lambda_{t+1}, T_{t+1}),
 \frac{4(t+2)}{\rho (\Lambda_{t+1}, \mathbb{Z}^{n+1} )}
  \right),
\label{newah}
\end{equation}
where $Z(\cdot ,\cdot ,\cdot, \cdot)$ is defined in
(\ref{starr}) and $H (\cdot ,\cdot ,\cdot, \cdot)$ is defined in (\ref{inte}).

Note that due to (\ref{ka})   we have
$$
 \mu_{k} (\xi_{t+1} , N) \le (2(t+1))^{-1},\,\,\,  H_t \le N \le H_t '
$$
$H_{t+1}$ may be  greater than $H_t'$. But it
follows from the inequality $ T_t \ge E_5 (t) $ and   the
definition (\ref{ahtprim}) of $H_t'$  that
$$
Z(2(t+1),\gamma_t^*,\Gamma_{t+1}, T_t)
  \le H_t'
.  $$  Lemma 6 implies   that  for
$$
N \ge Z(2(t+1),\gamma_t^*,\Gamma_{t+1}, T_t)
  $$
  we have
  $$
 \mu_{k} (\xi_{t+1} , N) \le (2(t+1))^{-1}.
$$
Hence
$$
 \mu_{k} (\xi_{t+1} , N) \le (2(t+1))^{-1},\,\,\,  H_t \le N \le H_{t+1}.
$$
 On the other hand,  it follows from (\ref{kaplusdva}) that
$$ \mu_{k+2} (\xi_{t+1} , N) \ge 2(t+1),\,\,\, N \ge H_t^n.
$$

For every $l$  the function $ \mu_l (\xi , N)$ is a continuous
function in $\xi$ and $N$. So, there exists $ \eta_{t+1}
>0
 $,
such that
\begin{equation}
 \mu_{k} (\xi , N) \le (t+1)^{-1},\,\,\,
\forall \xi : \, |\xi - \xi_{t+1} |\le \eta_{t+1}, \,\,\,\forall
N:\,
 H_t \le N \le H_{t+1},
\label{kaa}
\end{equation}
\begin{equation}
\mu_{k+2} (\xi , N) \ge t+1,\,\,\,   \forall \xi : \, |\xi -
\xi_{t+1} |\le \eta_{t+1},  \,\,\,\forall N:\, H_t^n \le N \le
H_{t+1}^n, \label{kaplusdva1}
\end{equation}
and
\begin{equation}
{\cal B}_{t+1} = \{ \xi = (1,\xi_1,...,\xi_n )\in \mathbb{R}^{n+1}
:\, |\xi - \xi_{t+1} |\le \eta_{t+1}\} \subset {\cal B}^*_t
\subset {\cal B}_{t} \label{ball1}
\end{equation}

 Now for
the objects $ \eta_i , \gamma_i,  W_i, H_i ,\Lambda_i ,
\Gamma_{i}, \xi_i $ with
 $i = t+1$ we have the following statements.

The condition
 (i) is satisfied by the construction.

The condition
 (ii) is satisfied
since $ \xi_{t+1}$ is a $(\Lambda_{t+1},\gamma_{t+1},
W_{t+1})$-BAD vector.

The condition
 (iii)
follows from (\ref{ball}).

The condition
 (iv) follows from (\ref{ball1}).

The condition
 (v) follows from  the definition (\ref{newah}) of $H_{t+1}$.

The condition
 (vi) follows from (\ref{ka}).

The condition
 (vii) follows from (\ref{kaplusdva1}).

The inductive procedure is described completely and Theorem 1 for
$ k \ge 2 $ is proved.

\vskip+1.0cm

The author thanks the anonymous referee for important suggestions.

\vskip+2.0cm

author: Nikolay Moshchevitin

\vskip+0.5cm

e-mail: moshchevitin@mech.math.msu.su, moshchevitin@rambler.ru

\end{document}